\documentclass[a4paper,11pt]{article}
\usepackage{graphicx}
\usepackage{pstricks}
\usepackage{amsmath,amssymb,amsfonts}
\usepackage{a4wide,epic}
\usepackage{color}

\newcommand{\lm}{\underline{\lambda}}
\newcommand{\lp}{\bar{\lambda}}

\begin{document}

\title{Adding a single state memory optimally accelerates symmetric linear maps}
\author{A.~Sarlette --- alain.sarlette@inria.fr\\ QUANTIC project team, Paris Sciences Lettres University,\\ INRIA Rocquencourt, Paris, France\\ \& SYSTeMS, Ghent University, Belgium.}
\date{June 2015}

\maketitle

%

\maketitle

\begin{abstract}
Previous papers have proposed to add memory registers to the dynamics of discrete-time linear systems in order to accelerate their convergence. In particular, it has been proved that adding one memory slot per agent allows faster convergence towards average consensus. We here prove that this situation \emph{cannot} be improved by adding more memory slots, when the knowledge about the self-adjoint linear map to be accelerated reduces to bounds on its extreme eigenvalues.
\end{abstract}



\section{Main result \& related work}\label{sec:probdef}


\subsection{Motivating question: consensus}

A consensus algorithm in discrete time \cite{TsitsiklisThesis} writes, for a column vector $x = [x_1;\,x_2;\,...; x_N] \in \mathbb{R}^N$:
\begin{eqnarray}\label{eq:base}
x(t+1) & = & x(t) - \alpha L \, x(t) \, .
\end{eqnarray}
Here $x_k \in \mathbb{R}$ is the state of agent $k$, $\alpha>0$ is a gain, and $L$ is the Laplacian matrix characterizing interactions among agents: component $j$ of $L \, x$ equals $\sum_{k=1}^N \, w_{j,k} \, (x_j-x_k)$ with weights $w_{j,k} \geq 0$. Usually $L$ contains many zeros, as each agent interacts with few fellows. For the algorithm to preserve the average of the initial values $x_k(0)$, we assume $w_{j,k} = w_{k,j}$ for all $j,k$. The Laplacian $L$ 
is then symmetric nonnegative definite, and if the interactions form a connected graph it has a single eigenvalue $\lambda_1 = 0$ with eigenvector $v_1 = [1;\, 1;\, ... 1]$; we call $v_1$ a consensus situation. When $(I_N-\alpha L)$ has nonnegative entries, with $I_N$ the $N\times N$ identity matrix, it can equivalently be viewed as the transition matrix of a Markov chain; symmetric $L$ implies that the limiting distribution is uniform.

For time-invariant $L$, the convergence of \eqref{eq:base} is dictated by the largest eigenvalue, in modulus, of $(I_N-\alpha L)$ --- excluding the trivial eigenvalue $1$ associated to $\lambda_1=0$ and eigenvector $v_1$, which spans the target subspace. In the orthonormal basis corresponding to the eigenvectors of $L$ (so-called ``modes''), the system decouples into
\begin{equation}\label{eq:basedec}
\tilde{x}_k(t+1) = (1-\alpha \lambda_k) \tilde{x}_k(t) \; , \;\; k=1,2,...,N,
\end{equation}
with $\tilde{x}_k$ the coefficient of mode $k$ and $\lambda_k$ the associated eigenvalue.
If we only know that the eigenvalues $\lambda_k$ of $L$ belong to an interval $[\lm,\lp] \subset \mathbb{R}_{>0}$ for $k=2,3,...,N$, then the convergence speed --- in terms of bounds on eigenvalues of \eqref{eq:basedec} guaranteed over all $\lambda_k \in [\lm,\lp]$ --- is optimal when $\alpha$ is selected to satisfy $(1-\alpha \lm) = -(1-\alpha \lp)$. This gives $\alpha = \frac{2}{\lp+\lm}$ and worst eigenvalue $\mu = \frac{\lp-\lm}{\lp+\lm}$.

In practice, the time step of discrete-time consensus is mostly limited by communication speed, much more rarely by local computation power. Hence \cite{Mea98} and later \cite{Joh08,Joh12,Cao13} propose to improve convergence speed by adding local dynamics to each agent, namely a memory slot with associated gain $\beta_1 \in \mathbb{R}$:
\begin{eqnarray}
\label{eq:mem1} x_k(t\!+\!1) = x_k(t) - \alpha {\textstyle \sum_{j=1}^N} \, w_{k,j} \, (x_k(t)-x_j(t))
 - \beta_1 (x_k(t) -x_k(t\text{-}1)) \, .
\end{eqnarray}
The papers analyze in detail only this case of one memory slot and prove that with $\beta_1 < 0$ it allows faster convergence. This spurs the natural question:\newline \textbf{\emph{\underline{How much can be gained by adding more memory slots?}}}\newline
The answer is~~\underline{\textbf{\emph{strictly nothing}}}~~, as we prove next in a more general setting.


\subsection{Formal result}\label{sec:result}

Consider a general linear iteration
\begin{equation}\label{eq:Ax}
x(t+1) = x(t) + \alpha (b - A x(t))
\end{equation}
where $x = [x_1;\,x_2;\,...; x_N] \in \mathbb{R}^N$ is a dynamic variable, $\alpha \in \mathbb{R}$ is a constant gain to be chosen, $b \in \mathbb{R}^N$ is a given constant bias and $A$ is a given constant, self-adjoint positive semidefinite linear map on $\mathbb{R}^N$. We are interested in applications where \eqref{eq:Ax} is applied iteratively in order to converge towards the set of fixed points $\mathcal{S} = \{ s \in \mathbb{R}^N : A s = b \}$. For consensus, with $A = L$ and $b=0$, we have $\mathcal{S} = \{ \lambda v_1 : \lambda \in \mathbb{R} \}$; in other computational applications $\mathcal{S}$ might reduce to a single point; we always assume $\mathcal{S}$ to be nonempty.

Motivated by the consensus setting, we consider a ``memory-accelerated'' version of \eqref{eq:Ax}:
$$ x(t+1) = x(t) + \alpha (b - A x(t)) + \sum_{m=0}^{M-1} \beta_m \, x(t\text{-}m) \, ,$$
with chosen gains $\beta_m \in \mathbb{R}$ for $m=0,1,2,...,M-1$. This yields stationary solutions $x$ satisfying $\alpha (b-Ax) + \left(\sum_{m=0}^{M-1} \beta_m\right) x = 0$. Hence we require $\sum_{m=0}^{M-1} \beta_m = 0$ to ensure that $\mathcal{S}$ remains a set of stationary solutions under the ``memory-accelerated'' dynamics. The latter can then be rewritten as:
\begin{eqnarray}\label{eq:AxPlus}
x(t+1) = x(t) + \alpha (b - A x(t)) + \sum_{m=1}^{M-1} \beta_m \, (x(t\text{-}m)-x(t)) \, ,
\end{eqnarray}
with freely chosen gains $\beta_m \in \mathbb{R}$. 
For analysis, we can rewrite \eqref{eq:AxPlus} in the basis of eigenvectors of $A$ and this just leads to a set of $N$ independent scalar equations
\begin{eqnarray}\label{eq:AxPlusB}
\tilde{x}_k(t+1) = (1-\alpha \lambda_k) \tilde{x}_k(t) \, + \alpha \tilde{b}_k  + \sum_{m=1}^{M-1} \beta_m \, (\tilde{x}_k(t\text{-}m)-\tilde{x}_k(t))
\end{eqnarray}
where $\lambda_k$ are the associated eigenvalues of $A$, for $k=1,2,...,N$. Hence it is clear that convergence of \eqref{eq:AxPlus} is governed by the separate convergence of each modal component $\tilde{x}_k$ associated to its eigenvalue $\lambda_k$ in \eqref{eq:AxPlusB}. Note that for $\mathcal{S}$ to be nonempty, we must have $\tilde{b}_k=0$ for every $k$ for which $\lambda_k=0$. Then the modes with $\lambda_k=0$ of both \eqref{eq:Ax} and \eqref{eq:AxPlus} (initialized with $\tilde{x}_k(t)=\tilde{x}_k(0)$ for $t<0$) become trivial, such that convergence towards $\mathcal{S}$ is governed by the \emph{nonzero} eigenvalues of $A$ only. We consider the setting where the choice of $\alpha,\beta_1,...,\beta_{M-1}$ is allowed to exploit the following knowledge. Let $\vert z \vert$ and $\angle(z)$ respectively denote the modulus and the phase of $z \in \mathbb{C}$.\vspace{2mm}

\noindent \textbf{Assumption 1:} \emph{It is known a priori that the nonzero eigenvalues $\lambda_k$ of $A$ take values in the interval $[\lm,\lp] \subset (0,+\infty)$.}\vspace{2mm}

\noindent \textbf{Definition 2:} \emph{Denote $\gamma_m(\lambda)$, $m=1,2,...,M$ the roots of the characteristic polynomial of the linear difference equation \eqref{eq:AxPlusB} with $\lambda_k = \lambda$. Then we define the \emph{convergence speed guarantee} $\nu$ of \eqref{eq:AxPlus} for fixed $\alpha,\beta_1,...,\beta_{M-1}$ by
$$\nu := \max_{\lambda \in [\lm,\lp]} \, \max_{m\in \{1,2,...,M\}} \, \vert \gamma_m(\lambda) \vert \, ,$$
i.e.~$\nu$ denotes the worst possible eigenvalue of \eqref{eq:AxPlusB} over all $\lambda_k$ satisfying Assumption 1.}\vspace{2mm}

We then prove the following result.\vspace{2mm}

\noindent \textbf{Theorem 3:} \emph{Consider dynamics \eqref{eq:AxPlus} with $A$ self-adjoint and satisfying Assumption 1. Then the convergence speed guarantee is optimized, i.e.~$\nu \in (0,1)$ is minimized over all $\alpha,\beta_1,\beta_2,...\beta_{M-1} \in \mathbb{R}$, when using just a single memory slot optimally tuned as in \cite{Mea98}, i.e.~taking:
\begin{eqnarray*}
\beta_m & = & 0 \, \text{ for all } m>1 \, ,\\
\beta_1 & = & \beta_{1*} := - \left( \frac{1}{\mu} - \sqrt{\frac{1}{\mu^2}-1} \right)^2 \, ,\\
\alpha & = & \alpha_* := \frac{2(1-\beta_{1*})}{\lp+\lm} \, , 
\end{eqnarray*}
where $\mu = \frac{\lp-\lm}{\lp+\lm}$. The corresponding convergence speed guarantee is
$$\nu_* = \frac{1}{\mu} - \sqrt{\frac{1}{\mu^2}-1} = \sqrt{-\beta_{1*}}\, .$$
}
This $\nu_*$ increases with $\mu$, which was the optimal convergence speed guarantee for $M=1$ i.e.~no added memory slot. For $\mu=1\text{-}\varepsilon$, with $\varepsilon$ frequently called the \emph{spectral gap}, we have $\nu_* < 1\text{-}\sqrt{\varepsilon}$ and the bound gets tight as $\varepsilon \rightarrow 1$.
\vspace{2mm}

\noindent \textbf{Remark 4:} \emph{\cite{Mea98} and related work mention that $\beta_1<0$ is necessary to get acceleration over the memoryless case. As an auxiliary result, it is not difficult to show for \eqref{eq:AxPlus} that if we are restricted to $\beta_m \geq 0$ for all $m$ then no acceleration with respect to the memoryless case can be achieved.}
\vspace{2mm}


\subsection{Related acceleration approaches}

Before proving Theorem 3 on the basis of complex polynomial analysis, we mention a few related acceleration schemes in the literature. We hope that in the future these viewpoints might lead to a more elegant explanation and intuitive understanding of the limitation to one memory slot.
 
\paragraph{Optimization techniques} The scheme studied here is not unlike techniques proposed earlier by the optimization community, where acceleration methods are a major focus point. For instance the Nesterov method \cite{Nesterov1983Method} accelerates the gradient descent optimization of a convex function $f$, i.e.~it accelerates 
$$x(t+1) = x(t) - \alpha \nabla_x f(x(t))$$
by applying, in its most basic version:
\begin{eqnarray*}
m(t) & = & x(t) - \alpha \nabla_x f(x(t))\\
x(t+1) & = & m(t) + g(t) \, (m(t)-m(t-1)) \, .
\end{eqnarray*}
For $f = \tfrac{1}{2}x^T A x - b^T x$, this corresponds exactly to the acceleration \eqref{eq:AxPlus} of $\eqref{eq:Ax}$ with $M=2$, except that the Nesterov method involves a time-varying step size $g(t)$, whose details go beyond the scope of the present paper. Our result seems to suggest that just adding more memory slots would not further accelerate the Nesterov method.

\paragraph{Robust control} The problem setting \eqref{eq:AxPlusB} with $\tilde{b}_k = 0$ can be viewed as an ensemble of closed-loop systems resulting from proportional feedback with different gains $\lambda_k$, i.e.:
\begin{eqnarray*}
y(z) & = & H(z) \, u(z) \; , \quad u(z) = - \lambda_k y(z)  \quad \text{with}\\[1mm]
H(z) & = & \frac{\alpha}{z-1 + {\textstyle \sum_{m=1}^{M-1}}\; \beta_m \, (1-z^{-m})}
\end{eqnarray*}
where $y(z)$ and $u(z)$ are the output and input respectively and $H(z)$ is the plant transfer function. We then want to design the plant $\alpha,\beta_1,...,\beta_{M-1}$ to get the fastest possible worst-case performance over the ensemble of feedback gains $\lambda_k \in [\lm,\lp]$. In view of Theorem 3, a plant with $M=2$ would be best in terms of convergence speed. As a major caveat though, the present paper proves optimality over all $H(z)$ with no zeros and with at least one pole at $z=1$; the possible practical relevance of this setting in control applications would have to be checked.

Various reformulations of where the uncertainty sits can of course be envisioned, e.g.~the uncertain plant could be modeled by $A$ whereas $\alpha,\beta_1,...,\beta_{M-1}$ would characterize the controller. From that perspective, \eqref{eq:AxPlusB} can be viewed as an optimal control problem for a so-called (non-symmetric) interval matrix or linear parametric uncertainty system, along the lines of Kharitonov's theorem. Unfortunately the latter does not generalize easily to discrete-time systems \cite{DTK,Vaid90}. An interesting result when there is a single uncertain parameter is proposed in \cite{Barmish}. Nevertheless, the author found no simple way to use this towards proving Theorem 3.

\paragraph{Accelerated consensus} As mentioned in the introduction, our initial motivation was the study of consensus algorithms accelerated through ``local'' memories; see \cite{Mea98,Joh08,Joh12,Cao13} which propose the scheme (similar to) \eqref{eq:AxPlus} and analyze the case $M=2$. Very recently \cite{Olshevsky2015} has proposed a scheme with ``extrapolating'' step, similar to the case $M=2$ (and to the Nesterov method, see above) but where the limited information known about the network is a bound on the total number of nodes, instead of on the extremal eigenvalues of $L$.

Researchers have also considered accelerations based on different resources than local memory, for instance allowing \emph{periodically time-varying} gain $\alpha$ in \eqref{eq:Ax}. Optimal tuning rules for $\alpha(t)$ can then build on \emph{polynomial-based filtering}~\cite{PolFt,SH07}. The optimal acceleration strategy of \cite{Mon13} can also be viewed as an instance of this framework, provided one limits the algorithm to a finite number of Chebyshev polynomials which are then applied periodically. In this latter context, when $\lambda_k$ can take all values in $[\lm,\lp]$, it is easy to check that the memory-based strategy \eqref{eq:AxPlusB} is faster than the optimal periodically time-varying one. However if more is known about the $\lambda_k$, then polynomial-based filtering is straightforward to adapt towards better speed-up. In the extreme case where $L$ has $M$ different eigenvalues whose values are exactly known, it is possible to construct a polynomial of order $M$ that achieves deadbeat convergence in $M$ steps. Similar ``finite-time consensus'' strategies have been proposed and studied in e.g.~\cite{Vkb13,Greek11}, with an interpretation of the resulting algorithms as gathering initial state information at a central node, computing its average and redistributing the average through the network.

It seems not straightforward to compare the resources ``local memory'' and ``time-varying gain'' outside an applicational context. In particular a standard time-varying implementation requires a stronger type of synchronization among all the nodes, namely they must not only agree on the frequency of updates but also on which gain to apply at each step. This may or may not be restrictive. The Nesterov method combines both memory and time-dependent step lengths. These approaches, along with accelerated Markov chains (see next paragraph), offer alternatives towards further acceleration but apparently \emph{only when more is known than $\lambda_k \in [\lm,\lp]$.}

\paragraph{Accelerated Markov chains} The optimal acceleration of the dual of consensus, namely Markov chains, has also been extensively studied in the literature, see e.g.~the ``lifting'' method of \cite{CLP99} and an even faster ``pseudo-lifting'' in \cite{JSS2010distributed}. Those approaches exploit knowledge of the particular network, hence they apply to the case where more is known than only $\lambda_k \in [\lm,\lp]$ (see previous paragraph). In fact \cite{JSS2010distributed} proposes an accelerated consensus version which directly connects with the  ``gather-and-distribute'' strategy of \cite{Vkb13,Greek11}, although with \emph{a priori} more complex communication requirements: it requires each node to communicate a vector of several values at each time step.

Yet lifted Markov chains cannot be viewed straightforwardly as the dual of accelerated consensus.\footnote{In fact, although this is a minor point, even in the unaccelerated case, Markov chains unlike consensus are restricted to systems with positive coefficients in the matrix $I-\alpha L$.} Whether the result of \cite{Mea98} and/or Theorem 3 can be used to investigate Markov chain acceleration, possibly without requiring detailed network knowledge, remains to be studied. The lifted Markov chain \cite{CLP99} \emph{using} detailed network knowledge is limited to quadratic speed-up.

\paragraph{Information theory} Consider the consensus application from the viewpoint of an individual node $k$. The messages received by node $k$ when \emph{symmetrically} exchanging data with a network, reflect in part the other nodes' initial values $\{ x_j(0) : j\neq k \}$ and in part the influence of $x_k(0)$ on the other nodes in the network. It seems not unreasonable to suspect that local memory can help disentangle these influences. Yet Theorem 3 says that if the network is poorly characterized, in the precise sense $\lambda \in [\lm,\lp]$, then having more than one additional memory is not helpful. This somehow seems to say: taking into account the direct feedback loop from $k$ through its neighbors $j$ and directly back to $k$ does allow improved convergence; but speculating about longer feedback loops does not pay off.\\

Finally, let us emphasize that Theorem 3 starts with a positive semidefinite self-adjoint matrix $A$, which can be viewed as a diffusive operator, or in consensus as an undirected communication matrix. In contrast, the state matrix characterizing the accelerated system \eqref{eq:AxPlus} is not symmetric. Hence it seems to introduce transport, or hidden directed communication (effects), which is known to improve convergence speed in several contexts \cite{CLP99,Barooah1}. This also implies that the acceleration of \cite{Mea98}, exactly like Markov chain liftings, cannot be applied iteratively. I.e.~the new state matrix, characterizing the system that results from an acceleration according to \cite{Mea98}, is not symmetric which precludes applying the technique of \cite{Mea98} once more.\\

In the following we first prove Theorem 3, before discussing an example in Section \ref{sec:example}.


\section{Proof, first part: reformulation}\label{sec:reform}

We consider the dynamics with memory slots in the form \eqref{eq:AxPlusB}. The convergence speed for each mode is governed by the roots of
\begin{eqnarray}\label{eq:nowold}
P(z;\lambda_k) = z^{M} - (1-\alpha \lambda_k) z^{M-1} + \; {\textstyle \sum_{m=0}^{M-2}} \; \beta_{M-m-1} \, (z^{M-1} - z^m) \, ,
\end{eqnarray}
viewed as a function of $z$ parameterized by $\lambda_k$.
In accordance with Definition 2, for the sequel we replace the set of $\lambda_k$ by a generic $\lambda$ which runs through $[\lm,\lp]$.


\subsection{Optimal solution with one memory slot}\label{ssec:opt2}

To make the paper self-contained, we outline a proof of optimal tuning for the case $M=2$, which is covered in \cite{Mea98}.\vspace{2mm}

\noindent \textbf{Proposition 5} [adapted from \cite{Mea98}]: \emph{For $M=2$ the tuning $\beta_1 = \beta_{1*}$, $\alpha = \alpha_*$, as proposed in Theorem 3, minimizes the value of the convergence speed guarantee $\nu$.
Moreover, the roots of 
$$P_*(z;\lambda) = z^2 - (1-\alpha_* \lambda) z + \beta_{1*}(z-1)$$
take the values $z_{\pm*}(\lambda) = \nu \, e^{\pm i \theta_{\lambda}}$ where the map $\lambda \mapsto\theta_{\lambda}$ is a continuous bijection from $[\lm,\lp]$ to $[0,\pi]$.}\vspace{2mm}

\noindent \underline{Proof:} First note that $\beta_1 \geq 1$ always leads to an unstable eigenvalue. Then the roots of $P(z;\lambda)$ are equal to the roots $z_{\pm}(\tilde{\lambda})$ of
$$f(z;\tilde{\lambda}) = z^2-(1-\beta_1) \tilde{\lambda} z - \beta_1$$
with $\tilde{\lambda} = (1- \alpha \lambda - \beta_1)/(1-\beta_1)$.  

We first fix a value of $\beta_1$ and study the roots of $f(z;\tilde{\lambda})$ as a function of $\tilde{\lambda}$, which translates into an optimal value for $\alpha$. Standard function analysis yields:
\newline \underline{Property 5:} \emph{The worst root $\max \{\vert z_+(\tilde{\lambda})\vert, \vert z_-(\tilde{\lambda})\vert \}$ of $f(z;\tilde{\lambda})$ is a monotonically (non-strictly) increasing function of $\vert \tilde{\lambda} \vert$.}
\newline Hence for given $\beta_1$ the choice of $\alpha$ should strive to minimize $\max_{\lambda \in [\lm,\lp]} (\vert \tilde{\lambda} \vert)$. This is obtained by choosing $\alpha$ such that $(1- \alpha \lm - \beta_1) = - (1- \alpha \lp - \beta_1)$, since $\vert \tilde{\lambda} \vert$ takes its extremal value(s) when $\lambda$ is at the boundary of the interval $[\lm,\lp]$. We hence get the condition $\alpha = \frac{2(1-\beta_{1})}{\lp+\lm}$. Moreover, this implies that $\tilde{\lambda}$ runs through the interval $[-\mu,\mu]$ when $\lambda$ spans $[\lm,\lp]$, with $\mu = \frac{\lp-\lm}{\lp+\lm}$ independent of $\beta_1$.

It remains to optimize $\beta_1$, assuming that it is associated to its respective optimal $\alpha$; by Property 5 and the previous sentence, this reduces to:
$$\min_{\beta_1 \in \mathbb{R}}\, \max_{\tilde{\lambda} \in \pm \mu} \max \left(\,\vert z_{\pm}(\tilde{\lambda}) \vert: f(z_{\pm}(\tilde{\lambda});\tilde{\lambda})=0 \,\right) \, .$$
A standard function analysis shows that it is optimal to take $\beta_1$ such that $z_+(\mu) = z_-(\mu) = - z_+(-\mu) = -z_-(-\mu)$. This yields the expressions of $\alpha_*, \beta_{1*}$ and $\nu$ given in Theorem 3. 

The just mentioned optimality condition requires that the discriminant $\Delta_* = (1-\beta_{1*})^2\tilde{\lambda}^2 + 4 \beta_{1*}$ of $f(z;\tilde{\lambda})=0$ takes its zeros for $\tilde{\lambda}=\pm \mu$. Then for all $\tilde{\lambda} \in [-\mu,\mu]$, i.e.~all $\lambda \in [\lm,\lp]$, we have $\Delta_* \leq 0$. Therefore $4 \vert z_{\pm*}(\lambda) \vert^2 = (1-\beta_{1*})^2 \tilde{\lambda}^2 -\Delta_*$ (sum of squared real and imaginary parts) $= -4\beta_1$ independently of $\lambda \in [\lm,\lp]$. Continuity of the roots of $f(z)$ between $\tilde{\lambda} = \mu$ and $\tilde{\lambda} = -\mu$ yields the last part of the Proposition.\hfill $\square$


\subsection{General polynomial property}

We now reformulate $P(z;\lambda)$ using Proposition 5.\vspace{2mm}

\noindent \textbf{Claim 6} [to be proved]: \emph{For any $\nu \in (0,1)$ and $\tilde{P}_{M-1} = \sum_{k=0}^{M-1} \, a_{k} y^k$ with $a_1,a_2,...,a_{M-1} \in \mathbb{R}$ and $a_{M-1} \neq -1$, there exists $\theta \in [0,\pi]$ for which
\begin{eqnarray}\label{eq:thepol}
\tilde{P}(y;\theta) = (y- e^{i\theta})(y- e^{-i\theta})\, y^{M-2} + \; (y-\tfrac{1}{\nu}) \; \tilde{P}_{M-1}(y)
\end{eqnarray}
has a root of modulus $\geq 1$.}\vspace{2mm}

\noindent \textbf{Proposition 7:} \emph{Theorem 3 is true if Claim 6 holds.}\vspace{2mm}

\noindent \underline{Proof} [of Proposition 7]: Let
\begin{eqnarray*}
\hat{P}(z;\theta) = (z-\nu e^{i\theta})(z-\nu e^{-i\theta}) z^{M-2} + (z-1) \; {\textstyle \sum_{m=0}^{M-1}} \; \tilde{\beta}_{m}\, z^m \, .
\end{eqnarray*}
Using $z^2 - (1-\alpha_* \lambda) z + \beta_{1*}(z-1) = (z-\nu e^{i\theta_\lambda})(z-\nu e^{-i\theta_\lambda})$ from Proposition 5, it is straightforward to check that $\hat{P}(z;\theta_\lambda) = \tfrac{\alpha_*}{\alpha} P(z;\lambda)$ provided 
\begin{eqnarray*}
\tilde{\beta}_k & = & \tfrac{\alpha_*}{\alpha} \; {\textstyle \sum_{m=0}^k} \, \beta_{M-m-1} \;\; \text{ for } k\leq M\text{-}3 , \\
\tilde{\beta}_{M\text{-}2} & = & \tfrac{\alpha_*}{\alpha}\; ({\textstyle \sum_{m=0}^{M-2}} \, \beta_{M-m-1}) \; - \beta_{1*} \, , \\
\tilde{\beta}_{M\text{-}1} & = & \tfrac{\alpha_*}{\alpha}-1 \, .
\end{eqnarray*}
Thus by taking appropriate $\tilde{\beta}_0,\tilde{\beta}_1,...\tilde{\beta}_{M-1} \in \mathbb{R}$, the roots of $\hat{P}(z;\theta_\lambda)$ can be made equal to the roots of $P(z;\lambda)$ from Theorem 3 for any $\alpha,\beta_1,...,\beta_{M-1} \in \mathbb{R}$, except for $\alpha=0$. The latter case is trivially uninteresting since for $\alpha=0$ any $x \in \mathbb{R}^N$ would be a stationary point of \eqref{eq:AxPlus}, 
hence no convergence can be obtained. 
Moreover, the case $\tilde{\beta}_{M-1} = -1$ will never appear since this would require infinite $\alpha$.
The tuning proposed in Theorem 3, with roots $\nu \, e^{\pm i \theta_\lambda}$ and $0$, corresponds to $\tilde{\beta}_0=\tilde{\beta}_1=...=\tilde{\beta}_{M-1}=0$. The parameter $\lambda \in [\lm,\lp]$ in $P(z;\lambda)$ is transferred to $\theta_\lambda \in [0,\pi]$ in $\hat{P}(z;\theta_\lambda)$, according to Proposition 5.

Hence a sufficient condition to prove Theorem 3 is to show that for any $\tilde{\beta}_0,\tilde{\beta}_1,...\tilde{\beta}_{M-1} \in \mathbb{R}$, with $\tilde{\beta}_{M-1} \neq -1$, there exists some $\theta \in [0,\pi]$ such that $\hat{P}(z;\theta)$ has at least one root with modulus $\geq \nu$. The statement of the Proposition is obtained by defining $y=z/\nu$ and $a_m = \tilde{\beta}_m \nu^{m-M+1}$ for $m=1,2,...,M-1$. \hfill $\square$\vspace{3mm}

Claim 6 looks like a discrete-time robust (in)stability property. Surprisingly, we know no standard result that would straightforwardly establish it. Examples can be constructed where the relevant roots are never real, or appear only for intermediate values of $\theta$ (see Section \ref{sec:example}), which seems to rule out simple variants of polynomial roots properties. Hence the remainder of this paper explicitly analyzes $\tilde{P}(y;\theta)$ in the complex plane. Note that for any polynomial $P(z)$, the modulus $\vert P(z) \vert$ is continuous for all $z \in \mathbb{C}$, while $e^{i \angle(P(z))}$ is continuous for all $z \in \mathbb{C} \setminus \{z : P(z) = 0 \}$.

\section{Proof, second part: analyzing $\tilde{P}(y;\theta)$}

We split up $\tilde{P}(y;\theta)$ into 
\begin{eqnarray*}
 P_1(y;\theta) & := & (y- e^{i\theta})(y- e^{-i\theta})\, y^{M-2} \; , \\
 P_2(y) & := & -(y-\tfrac{1}{\nu}) \, \tilde{P}_{M-1}(y) \; .
\end{eqnarray*}
Any $y$ at which $P_1(y;\theta) = P_2(y)$ is a root of $\tilde{P}(y;\theta)$. We first dispose of two special cases.\vspace{2mm}

\noindent \textbf{Proposition 8:} \newline \textbf{(a)} \emph{If $\tilde{P}_{M-1}$ in Claim 6 has a root on the unit circle, then $\tilde{P}(y;\theta)$ has a root $y_1$ with $\vert y_1 \vert \geq 1$ for some $\theta \in [0,\pi]$.}\newline
\textbf{(b)} \emph{If $a_{M-1}<-1$ in Claim 6, then $\tilde{P}(y;\theta)$ has a root $y_1$ with $\vert y_1 \vert \geq 1$ for all $\theta \in [0,\pi]$.}
\vspace{2mm}

\noindent \underline{Proof:} \textbf{(a)} Let $y_*$ the root of $\tilde{P}_{M-1}$ on the unit circle and take $e^{i\theta_*} = y_*$ if $\text{Imaginary}(y_*)\geq 0$, else $e^{-i\theta_*} = y_*$. Then $\tilde{P}(y_*,\theta_*) = 0$ with $\vert y_* \vert \geq 1$.

\textbf{(b)} We separately analyze (i) the sign and (ii) the magnitude of both $P_1$ and $P_2$ for $y$ belonging to the real positive axis.
\newline (i) Let $y_* \geq 1/\nu > 1$ the largest real root of $P_2(y)$. Both $P_2(y)$ and $P_1(y;\theta)$ are positive for $y \in (y_*,+\infty)$, for any $\theta$.
\newline (ii) $\vert P_2(y) \vert < \vert P_1(y;\theta) \vert$ for $y$ close to $y_*$, while for very large $\vert y \vert$ we have $\vert P_2(y) \vert \simeq \vert a_{M-1}\vert\, \vert y \vert^M > \vert P_1(y;\theta) \vert \simeq \vert y \vert^M$.
\newline Hence for any $\theta$, there exists some $y_1 \in (y_*,+\infty)$ where $\vert P_1(y_1;\theta) \vert = \vert P_2(y_1)\vert$, which with (i) implies $P_1(y_1;\theta) = P_2(y_1)$. This $y_1$ is a root of $\tilde{P}(y;\theta)$ with $\vert y_1 \vert > 1$.\hfill $\square$\\


From Prop.~8(a), we can reduce our investigation to the case where the roots of $P_2$ are disjoint from the roots of $P_1$. Indeed, we have excluded roots of $P_2$ on the unit circle, whereas a situation with $P_2$ and $P_1$ having the common root $y=0$ can be reduced to an equivalent algorithm with a lower value of $M$.

The analysis of $\tilde{P}(y;\theta)$ towards Claim 6 for the remaining cases is essentially a generalization of the proof of Prop.~8(b), from the real line towards the complex plane.


\subsection{Partitioning the complex plane with $P_1$ and $P_2$}\label{ssec:construction}

For given $P_1(\cdot;\theta)$ and $P_2(\cdot)$, we can partition the complex plane into a collection of open connected sets where $\vert P_1 \vert > \vert P_2 \vert$ (we call these type 1 sets), a collection of open connected sets where $\vert P_2 \vert > \vert P_1 \vert$ (we call these type 2 sets), separated by sets where $\vert P_2 \vert = \vert P_1 \vert$. Note that, as the roots of $P_1$ and $P_2$ can be assumed disjoint (see Prop.8(a)), any root of $P_2$ unambiguously belongs to a type 1 set and any root of $P_1$ belongs to a type 2 set.

\begin{figure*}
\includegraphics[trim=4cm 11cm 4cm 0.5cm, clip=true, width=110mm]{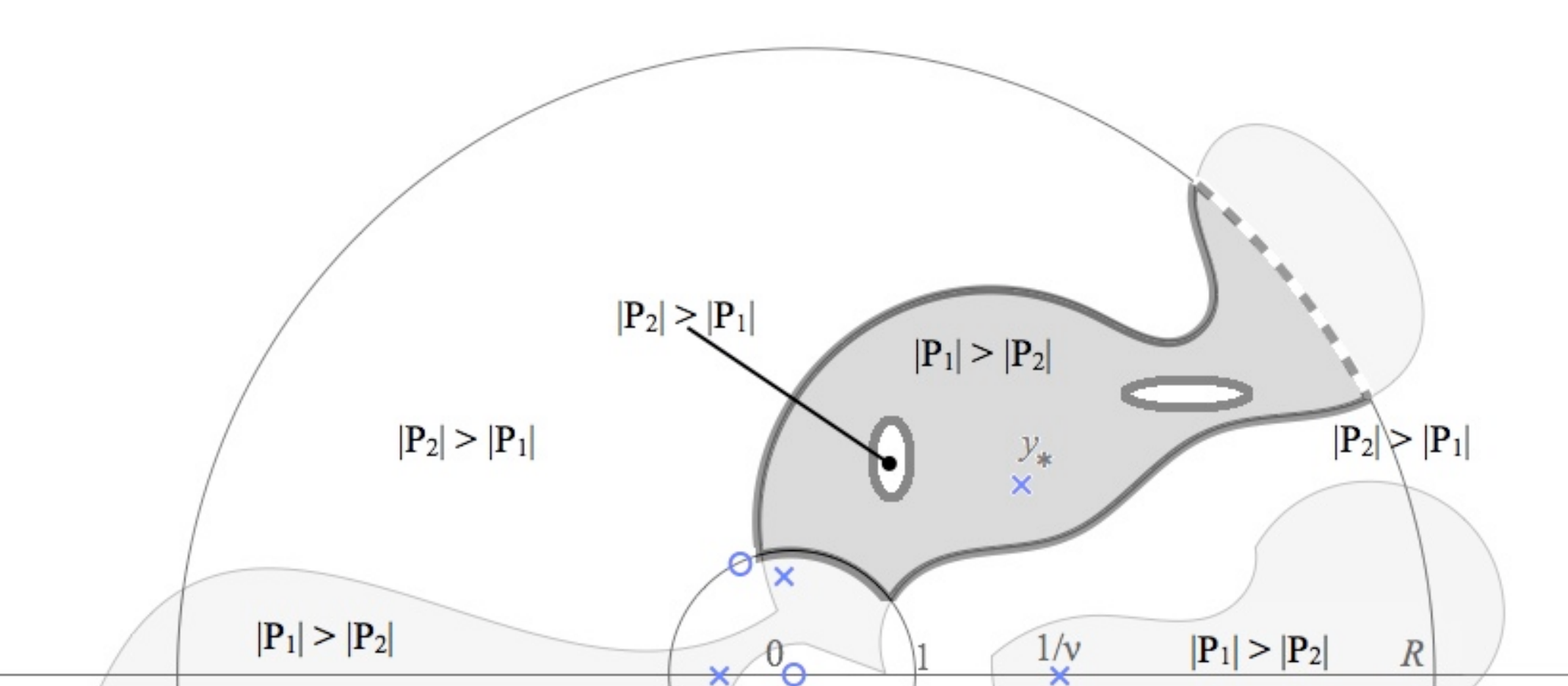}\\
\caption{Illustrating the definition of the type 1 sets (shaded), type 2 sets (plain white) and $\Gamma_R(y_*)$ (darker shading). Small blue crosses represent roots of $P_2$, small blue circles the roots of $P_1$. The boundary $\partial \Gamma_R(y_*)$ (thick plain and dotted curves) and its subset $\partial\Gamma_{\bar R}(y_*)$ (thick plain curves) are also represented. The complete picture is meant to be symmetric around the real axis. The sets drawn here are purely schematic, they are not meant to correspond to any actual polynomials $P_1$, $P_2$ and may possibly be more general than what polynomial properties would allow.}\label{fig:construction}
\end{figure*}

For any root $y_*$ of $P_2$ with $\vert y_* \vert > 1$, let $\Gamma(y_*)$ the type 1 set containing $y_*$. Then for any $R>\vert y_* \vert$, let
$$
\tilde{\Gamma}_R(y_*) := \Gamma(y_*) \cap \{y \in \mathbb{C} : 1 < \vert y \vert < R \}
$$
and define $\Gamma_R(y_*)$ to be the connected component of $\tilde{\Gamma}_R(y_*)$ such that $y_* \in \Gamma_R(y*)$ (see Figure \ref{fig:construction}). There is at least one such set corresponding to $y_* = 1/\nu$. We denote the boundary of $\Gamma_R(y_*)$ by $\partial \Gamma_R(y_*)$. Note that for any $y \in \partial \Gamma_R(y_*)$ we have either $\vert P_2(y) \vert = \vert P_1(y;\theta) \vert$, or $\vert P_2(y) \vert < \vert P_1(y;\theta) \vert$; the latter case involves points where either $\vert y \vert = 1$, or $\vert y \vert = R$. We finally define $\partial\Gamma_{\bar R}(y_*) = \partial \Gamma_R(y_*) \cap \{ y : \vert y \vert < R \}$ and reduce the proof of Claim 6 to the following.
\vspace{2mm}

\noindent \textbf{Claim 9} [to be proved]: \emph{For any $P_2$ --- excluding the cases already handled by Proposition 8 --- there exists $\tilde{\theta} \in [0,\pi]$ and $y_*$ a root of $P_2$ such that $e^{i\angle(P_2(\tilde{y}))} = e^{i\angle(P_1(\tilde{y};\tilde{\theta}))}$ for some $\tilde{y} \in \partial \Gamma_{\bar R}(y_*)$.}
\vspace{2mm}

\noindent \textbf{Proposition 10:} \emph{Claim 6 is true if Claim 9 holds.}\vspace{2mm}

\noindent \underline{Proof:} We must show that for any situation satisfying Claim 9, we can find $y_1 \in \mathbb{C}$ with $\vert y_1 \vert \geq 1$ and $\theta \in [0,\pi]$, such that $P_1(y_1;\theta) = P_2(y_1)$. For given $P_2$, excluding the cases covered by Proposition 8, take $\tilde{y},\tilde{\theta}$ according to Claim 9.


If $\tilde{y}$ belongs to the part of $\partial \Gamma_{\bar R}(y_*)$ where $\vert P_2(y) \vert = \vert P_1(y;\tilde{\theta}) \vert$, we have $P_1(\tilde{y};\tilde{\theta}) = P_2(\tilde{y})$ with $\vert \tilde{y} \vert \geq 1$ by construction. Hence the case is closed with $y_1 = \tilde{y}$, $\theta=\tilde\theta$.

If $\vert P_2(\tilde{y}) \vert \neq \vert P_1(\tilde{y};\tilde{\theta}) \vert$, then necessarily $\vert P_2(\tilde{y}) \vert < \vert P_1(\tilde{y};\tilde{\theta}) \vert$ and $\vert \tilde{y} \vert = 1$. Let $\tilde{y}=e^{i\phi}$ and assume $\phi \in [0,\pi]$; the case $\phi \in [-\pi,0]$ is the same modulo a few notation changes. As for Prop.~8(b), we separately investigate (i) the phase and (ii) the modulus of $P_1$, $P_2$.
\newline (i) By computing
$$
(\tilde{y}-e^{i\theta})(\tilde{y}-e^{-i\theta})\tilde{y}^{M-2}\; e^{-i (M-1) \phi} = 2 \left( \cos(\phi) - \cos(\theta) \right).
$$
we see that $e^{i\angle(P_1(\tilde{y};\theta))} = e^{i\angle(P_1(\tilde{y};\tilde{\theta}))} = e^{i\angle(P_2(\tilde{y}))}$ for all $\theta$ for which $\left( \cos(\phi) - \cos(\theta) \right)$ takes the same sign as $\left( \cos(\phi) - \cos(\tilde{\theta}) \right)$, i.e.~for all $\theta \in [\tilde{\theta},\phi) =:\mathcal{I}_1$ or all $\theta \in (\phi,\tilde{\theta}]=:\mathcal{I}_1$, depending on the ordering of $\tilde{\theta}$ and $\phi$.
\newline (ii) $f(\theta) := \vert P_2(\tilde{y}) \vert - \vert P_1(\tilde{y};\theta) \vert$ is a continuous function of $\theta$, with $f(\tilde{\theta}) < 0$ and $f(\phi) > 0$ by construction. Hence there exists $\theta_1 \in \mathcal{I}_1$ such that $f(\theta_1) = 0$.
\newline Combining (i) and (ii), we have $P_1(y_1;\theta) = P_2(y_1)$ with $y_1 = \tilde{y}$ and $\theta=\theta_1$.\hfill $\square$


\subsection{Proving Claim 9}\label{ssec:Claim 9}

It remains to prove that Claim 9 is true, which in fact is a consequence of \emph{Cauchy's argument principle} in complex analysis. Adapted to the current case, the argument principle states the following fact.\vspace{2mm}

\noindent \textbf{Property 11} [Cauchy's Argument Principle]: \emph{Let $\mathcal{D} \subset \mathbb{C}$ a bounded, simply connected open set whose boundary $\partial\mathcal{D}$ is the image of a closed simple curve $\gamma: [0,1] \rightarrow \mathbb{C}$. Assume that $\mathcal{D}$ contains $n$ roots of $P_2(y)$ and $p$ roots of $P_1(y;\theta)$ for some fixed $\theta$, while $\partial \mathcal{D}$ contains no roots.
Then there exists a continuous function $f(t): [0,1] \rightarrow \mathbb{R}$ such that 
$$\mathrm{exp}\left(i \left( \vphantom{\int} \angle(P_1(\gamma(t);\theta))-\angle(P_2(\gamma(t))) \right)\right) = \mathrm{exp}(i f(t)) \; ,$$
and any such $f(t)$ satisfies $\vert f(1)-f(0) \vert = 2\pi \vert p-n \vert$.}\vspace{2mm}

We will apply this principle with $\partial\mathcal{D}$ defined by elements of $\partial\Gamma_{R}(y_*)$. Note that the boundary of any $\Gamma_R(y_*)$ is indeed sufficiently regular since the locus in the complex plane where $\vert P_1 \vert^2 - \vert P_2 \vert^2 = 0$ is by definition a planar algebraic curve. We will also use the following related property.\vspace{2mm}

\noindent \textbf{Property 12:} \emph{Consider the setting of Property 11 with $\vert p-n \vert =: m > 0$. Then there are at least $m$ points $y=\gamma(t)$ on $\partial\mathcal{D}$ where $e^{i\angle(P_2(y))} = e^{i\angle(P_1(y;\theta))}$.}\vspace{2mm}

\noindent \underline{Proof:} The interval $(f(0),f(0)+2\pi m] \subset \mathbb{R}$ contains $m$ different multiples of $2\pi$ and the function $f(t)$ must take each of these values at least once as it continuously evolves from $f(0)$ to $f(1) = f(0) + 2\pi m$. The $\gamma(t)$ associated to these values of $f(t)$ satisfy $e^{i\angle(P_2(y))} = e^{i\angle(P_1(y;\theta))}$.\hfill $\square$\\

\begin{figure*}
\hspace{0.4cm}
\includegraphics[trim=13cm 11cm 4cm 0.5cm, clip=true, width=75mm]{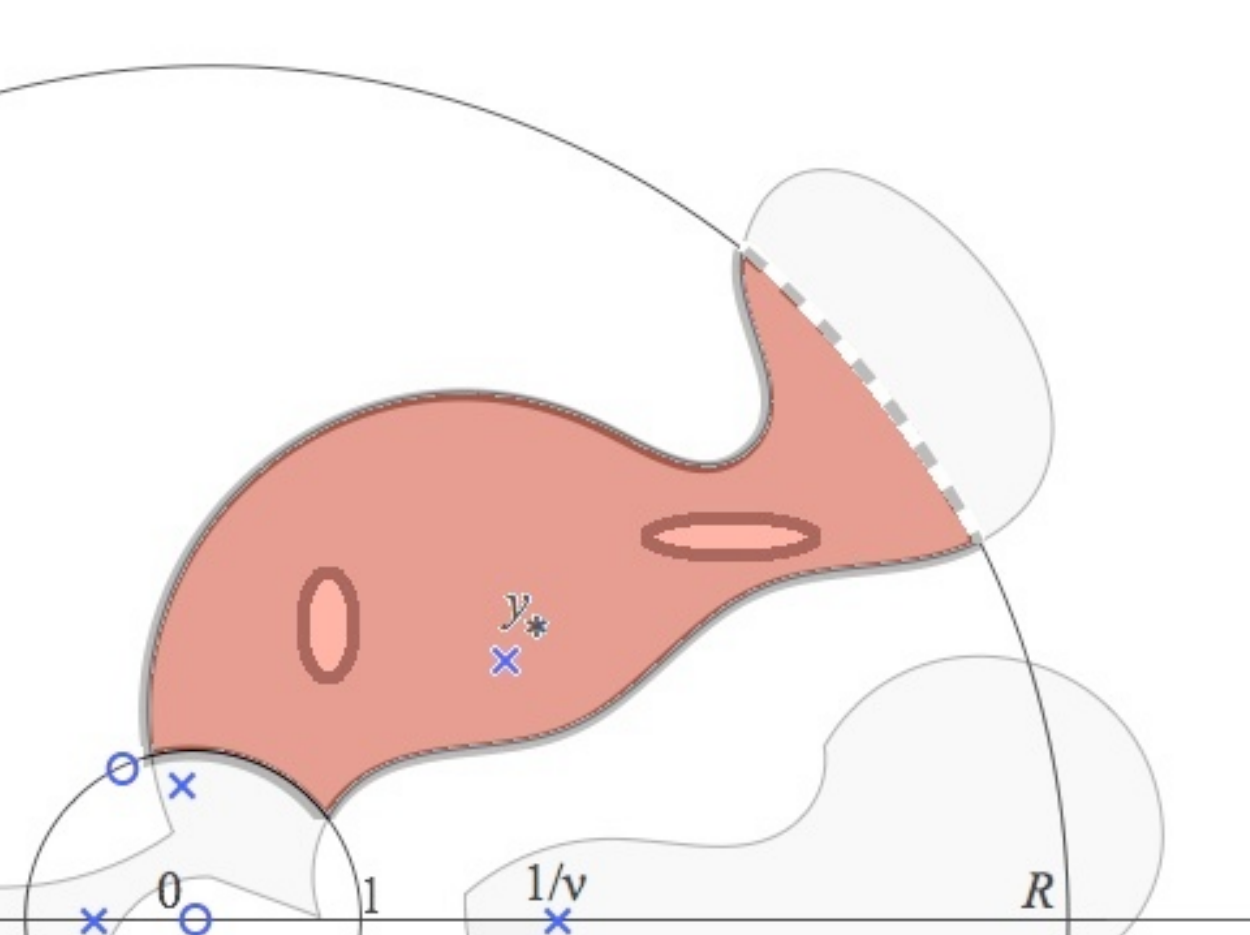} 
\includegraphics[width=80mm]{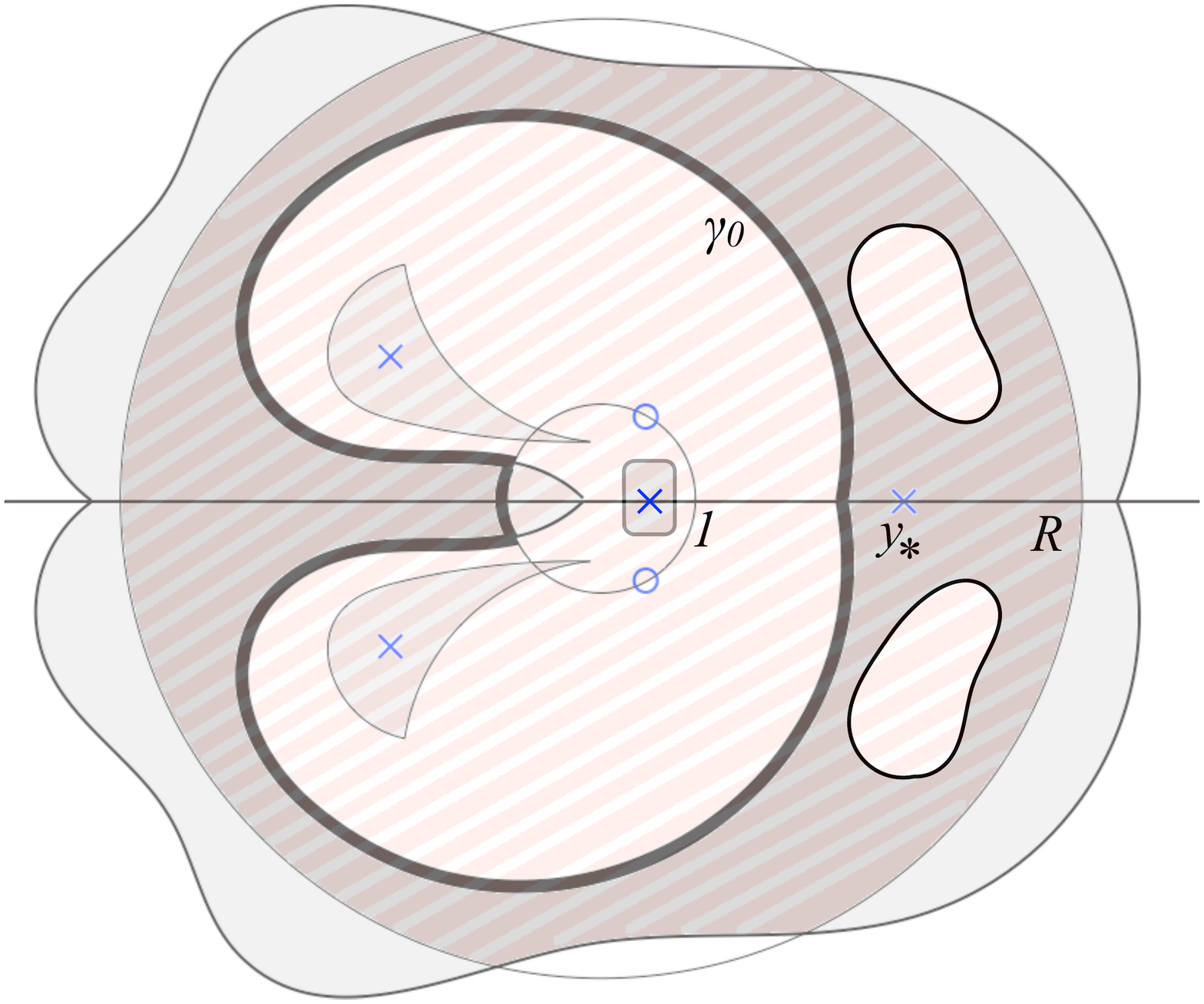}
\caption{Schematic illustration of Case A (left) and Case B (right) to which we apply Properties 11 and 12. The sets and roots are represented similarly to Fig.~\ref{fig:construction}. For Case A the set in red is $\mathcal{D}_R(y_*)$. For Case B, the striped set is $\mathcal{D}_R(y_*)$, the thick curve is $\gamma_0$ and the set enclosed by $\gamma_0$ is $\mathcal{D}_0(y_*)$.}\label{fig:Biggy2}
\end{figure*}

The remainder of the proof relies on the following observations (see Fig.~\ref{fig:Biggy2}). For some fixed $\theta$, let $\mathcal{D}_R(y_*)$ the smallest simply connected set containing $\Gamma_R(y_*)$, where the latter contains $m \geq 1$ roots of $P_2$. Note that $\partial\mathcal{D}_R(y_*) \subseteq \partial\Gamma_R(y_*)$.\vspace{1mm}
\newline \textbf{Case A:} If $\mathcal{D}_R(y_*)$ does not contain the open unit disc (see Fig.~\ref{fig:Biggy2} left), then it contains at least $m$ roots of $P_2$ but no root of $P_1$. Hence applying Properties 11 and 12 with $\mathcal{C} = \mathcal{D}_R(y_*)$, there must be $m$ points on $\partial\mathcal{D}_R(y_*) \subseteq \partial\Gamma_R(y_*)$ where $e^{i\angle(P_2(y))} = e^{i\angle(P_1(y;\theta))}$.
\newline \textbf{Case B:} If $\mathcal{D}_R(y_*)$ contains the open unit disc, then $\partial\Gamma_{\bar R}(y_*)$ contains a closed curve $\gamma_0$ that separates $\Gamma_R(y_*)$ from the open unit disc and we can define $\mathcal{D}_0(y_*) \subset (\mathcal{D}_R(y_*) \setminus \Gamma_R(y_*))$ the simply connected set whose boundary is $\gamma_0$ (see Fig.~\ref{fig:Biggy2} right). This set contains the $M$ roots of $P_1$ but at most $M-m$ roots of $P_2$. Hence by Properties 11 and 12 there must be $m$ points on $\gamma_0 \subset \partial\Gamma_{\bar R}(y_*)$ where $e^{i\angle(P_2(y))} = e^{i\angle(P_1(y;\theta))}$.
\vspace{1mm}\newline
Claim 9 is readily true for Case B. For Case A, we need to ensure that at least one of the $m$ points in $\partial\Gamma_R(y_*)$ also belongs to $\partial\Gamma_{\bar R}(y_*)$. Towards this, we investigate $e^{i\angle(P_1(y;\theta))-i\angle(P_2(y))}$ when $\vert y \vert =R$ with $R$ very large.\vspace{2mm}

\noindent \textbf{Lemma 13:} \emph{If in $P_2$ the coefficient $a_{M-1}>0$ then there exists $R_1>0$ such that\newline $e^{i\angle(P_1(y;\theta))-i\angle(P_2(y))} \neq 1$ for all $y$ with $\vert y \vert > R_1$ and for all $\theta \in [0,\pi]$.}\vspace{2mm}

\noindent \underline{Proof:}  Note that $e^{i\angle(P_1(y;\theta))-i\angle(P_2(y))} = \mathrm{exp}(i\angle\left( \frac{P_1(y;\theta)}{P_2(y)} \right))$. 
Let us evaluate the phase of the ratio of two general polynomials $P_g(y) = \sum_{k=0}^M g_k y^k$ and $P_h(y) = \sum_{k=0}^M h_k y^k$ with real coefficients $g_k, h_k$ and assuming $g_M h_M < 0$. Applying the formula $\frac{a+b i}{c+ d i} = \frac{1}{c^2+d^2}((ac+bd) + (bc-ad) i)$, with $a,b,c,d$ real and $i$ the imaginary unit, we get with $y=R e^{i\phi}$:
\begin{eqnarray*}
\frac{P_g(R e^{i \phi})}{P_h(R e^{i \phi})} = \frac{\sum_{k=0}^M \tfrac{g_{M-k}}{R^{k}} (\cos(k\phi) - i\sin(k\phi))}{\sum_{k=0}^M \tfrac{h_{M-k}}{R^{k}} (\cos(k\phi) - i\sin(k\phi))}
= \eta \; \left(g_M h_M + O(1/R) + i O(1/R) \right)
\end{eqnarray*}
for some $\eta>0$. Thus there exists $R_1>0$ such that 
$$\text{Real}\left(\tfrac{P_g(R e^{i \phi})}{P_h(R e^{i \phi})} \right) = g_M h_M + O(1/R) < 0$$ 
for all $R>R_1$.
The Lemma follows by applying the above result with $P_g(y)=P_1(y;\theta)$ and $P_h(y)=P_2(y)$.\hfill $\square$\\

We are now ready to prove Claim 9 in a few steps.\vspace{2mm}

\noindent \underline{Proof of Claim 9 for $a_{M-1} \leq 0$:} From Prop.8(c) we only need to consider the case $a_{M-1} \in (-1,0]$. Then there exists $R_1>1$ such that $\vert y \vert^M \simeq \vert P_1(y;\theta) \vert > \vert P_2(y) \vert \simeq \vert a_{M-1} \vert \, \vert y \vert^M$ for all $y$ for which $\vert y \vert > R_1$. Take any root $y_*$ of $P_2$ with $\vert y_* \vert >1$. Then for $R>R_1$, the points $y$ with $R_1 < \vert y \vert < R$ cannot belong to $\partial\Gamma_R(y_*)$, since this would require $\vert P_1(y;\theta) \vert = \vert P_2(y) \vert$. This implies that $\partial\Gamma_R(y_*)$ contains either all points or no point of the circle $\mathcal{C}_R := \{y\in\mathbb{C}:\vert y \vert = R\}$. In case $\partial\Gamma_R(y_*) \cap \mathcal{C}_R = \emptyset$, we have $\partial\Gamma_R(y_*) = \partial\Gamma_{\bar R}(y_*)$ and the observations after Proposition 12 allow to directly conclude. In case $\partial\Gamma_R(y_*) = \mathcal{C}_R$, we are necessarily in Case B of the observations after Proposition 12 hence the conclusion is also immediate.\hfill $\square$\\

\noindent \underline{Proof of Claim 9 for $a_{M-1}>0$:} Take $R>R_1$ to satisfy Lemma 13. The set $\partial\Gamma_{R}(y_*) \setminus \partial\Gamma_{\bar R}(y_*)$ is a subset of $\{ y \in \mathbb{C} : \vert y \vert = R \}$, and Lemma 13 implies $e^{i\angle(P_2(y))} \neq e^{i\angle(P_1(y;\theta))}$ for all $y$ for which $\vert y \vert = R$. Hence all the points identified in the observations after Proposition 12, where $e^{i\angle(P_2(y))} = e^{i\angle(P_1(y;\theta))}$ on $\Gamma_{R}(y_*)$, must belong to $\Gamma_{\bar R}(y_*)$.\hfill $\square$\\

With this we have covered all situations and hence concluded the proof of Claim 9, which proves Theorem 3.


\section{Example}\label{sec:example}

Consider a linear map $A$ with nonzero eigenvalues $\lambda_k \in [0.0122,0.9878]$. The optimal acceleration without added memory slot i.e.~with $M=1$ yields $\mu = 0.9756$ and a spectral gap $1-\mu = 0.0244$. With $M-1=1$ added memory slot, an improved convergence speed guarantee $\nu = 0.8000$ is obtained, with $1-\nu = 0.2 < \sqrt{0.0244} = 0.1562$, using the optimal parameters $\alpha_* = 3.2800$ and $\beta_{1*} = -0.6400$.\\

If more is known about the eigenvalues of $A$ then additional memory does help accelerate. For instance consider $A$ with nonzero eigenvalues $\lambda_k \in \Lambda = [0.0122,0.0182] \cup \{0.9878\}$ i.e.~the largest eigenvalue is in fact isolated far away from the others. This does not allow to improve convergence speed with $M=2$, because for $M=2$ the bound $\nu \geq 0.8000$ holds as soon as $A$ features both eigenvalues $0.0122$ and 0.9878 (see the proof of Prop.~5).

\begin{figure}
\includegraphics[trim=1cm 2cm 1cm 2cm, clip=true, width=90mm]{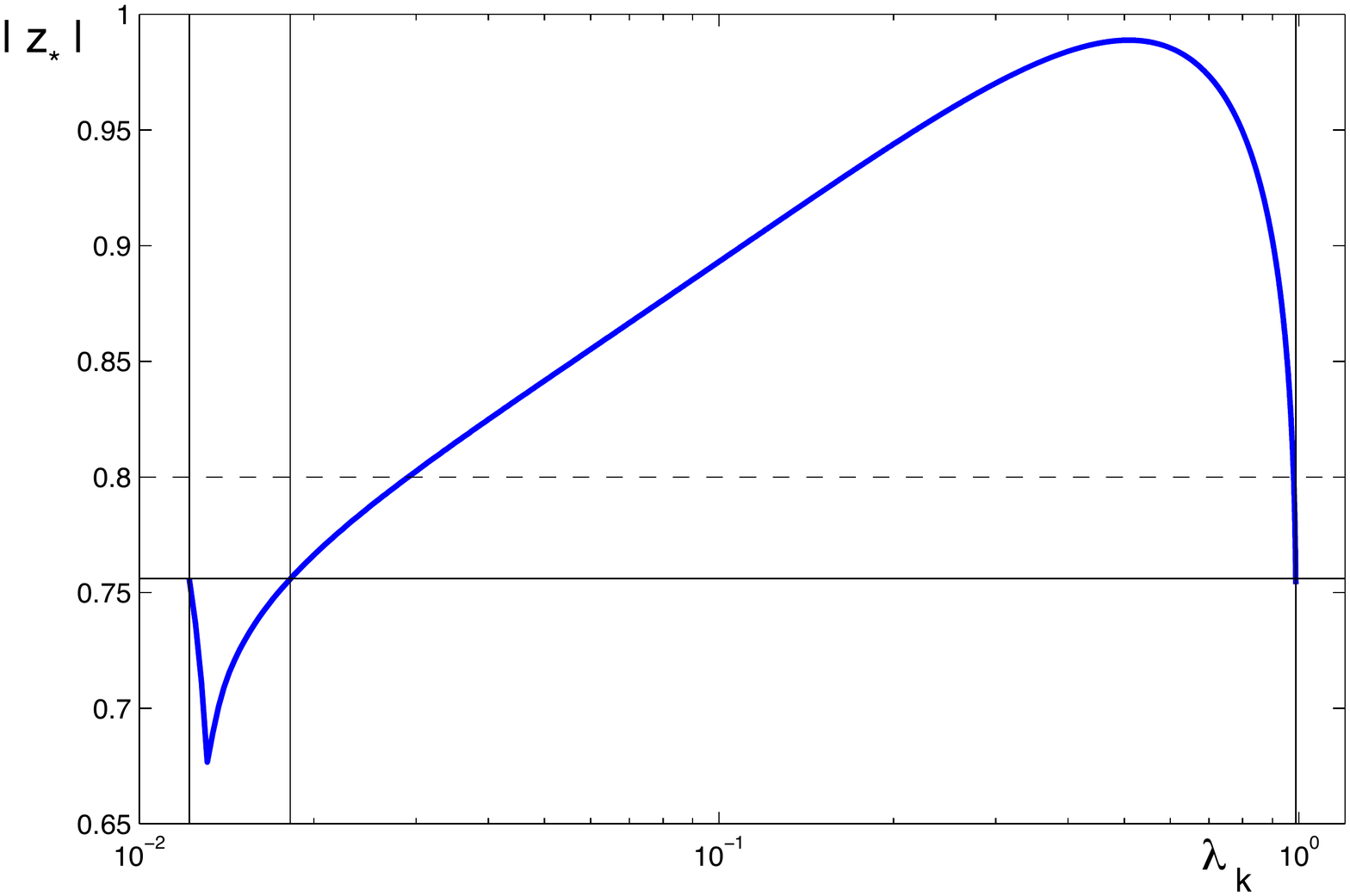}
\caption{Largest root in modulus $\vert z_* \vert$ of \eqref{eq:nowold} with $M=4$, $\alpha = 3.6908$, $\beta_1 = -0.9083$, $\beta_2=0.006662$, $\beta_3=0.06785$, as a function of $\lambda_k \in [0.0122,0.9878]$. Reducing the relevant eigenvalue set to $\lambda_k \in [0.0122,0.0182] \cup \{0.9878\}$ (vertical lines) yields an improved convergence speed guarantee $\tilde{\nu} = 0.7560$ (horizontal full line) with respect to the optimal situation $\nu = 0.8000$ with $M=2$ (horizontal dashed line). However this happens to the detriment of $\vert z_* \vert > 0.8000$ for all $\lambda_k \in [0.0291,0.9788]$.}\label{fig:lint}
\end{figure}

In contrast, $M=4$ and parameter values $\alpha = 3.6908$, $\beta_1 = -0.9083$, $\beta_2=0.006662$, $\beta_3=0.06785$ do yield an improved convergence speed guarantee $\tilde{\nu} = 0.7560$ over all $\lambda_k \in \Lambda$ (with $\tilde{\nu}$ defined by straightforward extension of Def.~2). Note that these parameters have been obtained numerically by local search around $\alpha = \alpha_*$, $\beta_1 = \beta_{1*}$ and $\beta_k=0$ for $k>1$; we do not exclude the existence of better parameter values. Figure \ref{fig:lint} shows $\vert z_* \vert$ the corresponding largest root in modulus of \eqref{eq:nowold} as a function of $\lambda_k$. It highlights how the improved $\vert z_* \vert$ for $\lambda_k \in \Lambda$ comes to the detriment of (much) worse $\vert z_* \vert$ for $\lambda_k \in [0.0291,0.9788]$.\\

\begin{figure}
\includegraphics[trim=12cm 7cm 10cm 7cm, clip=true, width=48mm]{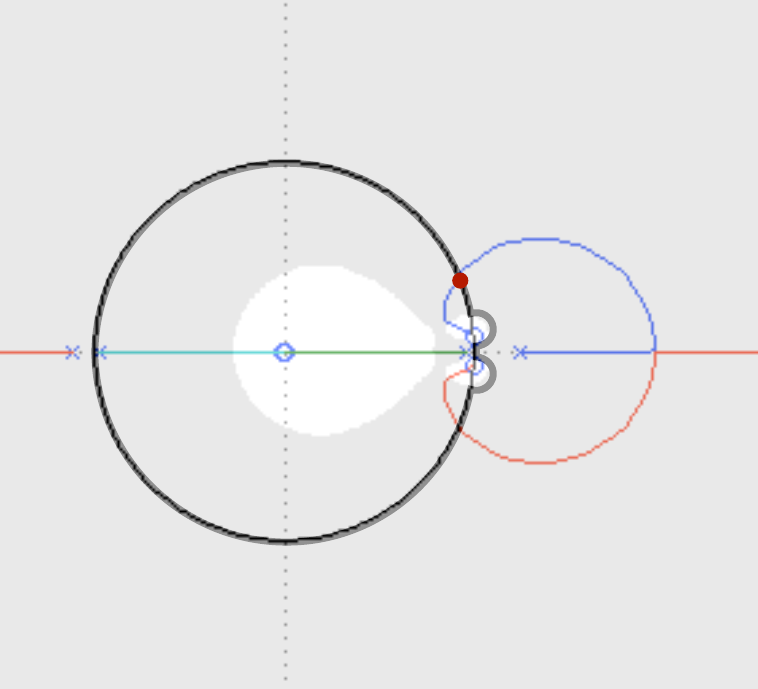} \hspace{2mm}
\includegraphics[trim=12cm 7cm 10cm 7cm, clip=true, width=48mm]{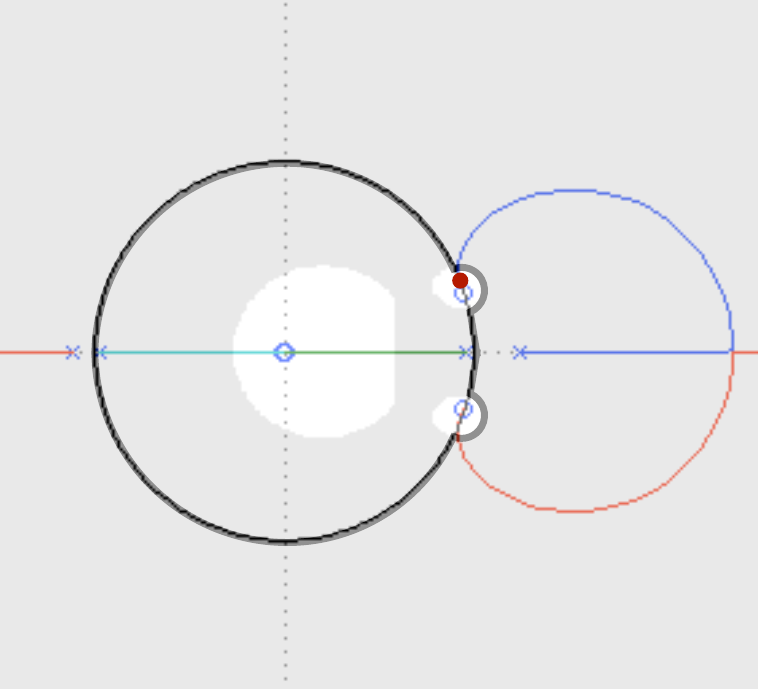} \hspace{2mm}
\includegraphics[trim=12cm 7cm 10cm 7cm, clip=true, width=48mm]{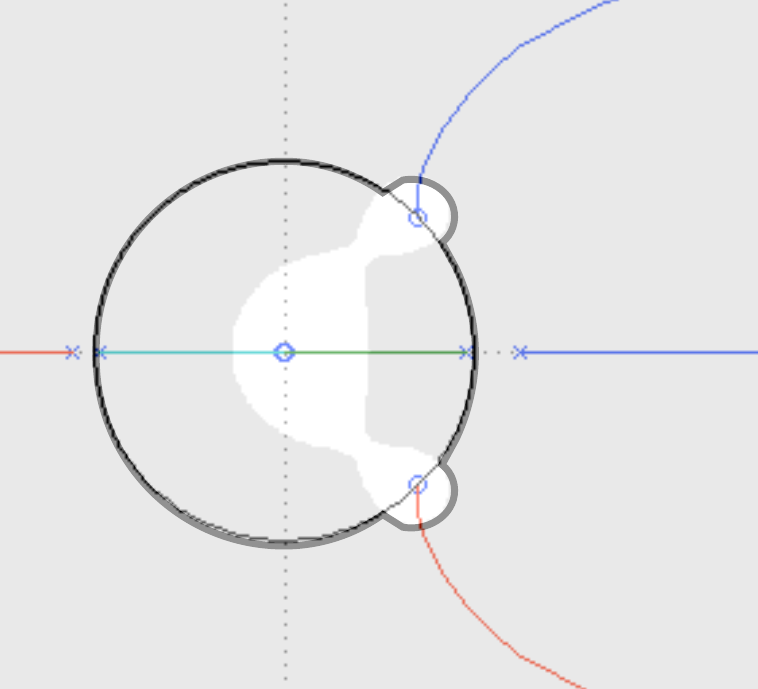}
\caption{Characteristics of $P_1(y;\theta)$ and $P_2(y)$ for three values of $\theta \in \{\, \pi/40,\,\pi/10, \, \pi/4 \, \}$, illustrating the proof of Claim 9 and Prop.10 on our example polynomial. See main text for details.}\label{fig:Ex}
\end{figure}

Finally, let us illustrate the proof of Claim 9 and Prop.10 with those values, i.e.~$M=4$ and $\alpha = 3.6908$, $\beta_1 = -0.9083$, $\beta_2=0.006662$, $\beta_3=0.06785$ for $[\lm,\lp] = [0.0122,0.9878]$ and hence $\nu=0.8000$. 
For three values of $\theta \in \{\, \pi/40,\,\pi/10, \, \pi/4 \, \}$, Figure \ref{fig:Ex} shows in white the set in $\mathbb{C}$ where $\vert P_1(y;\theta)\vert < \vert P_2(y) \vert$ and in grey the set where $\vert P_1(y;\theta)\vert > \vert P_2(y) \vert$. Small circles are the roots of $P_1$ and small crosses the roots of $P_2$. In the middle is the unit circle and the colored lines denote the locus where $e^{i\angle(P_1)} = e^{i\angle(P_2)}$. For all values of $\theta$, we are in Case B of Section \ref{ssec:Claim 9}. Correspondingly, we have highlighted in black the curve $\gamma_0$ delimiting the simply connected set $\mathcal{D}_0(y_*)$. For all $\theta$, there are two points along $\gamma_0$ where $e^{i\angle(P_1)} = e^{i\angle(P_2)}$. For $\theta=\pi/40$ the crossing occurs at some $y_1$ on the unit circle (red dot on Fig.~\ref{fig:Ex}), in a region where $\vert P_1 \vert > \vert P_2 \vert$. For $e^{i\theta}$ closer to $y_1$, e.g.~$\theta = \pi/10$, we still have $e^{i\angle(P_1(y_1;\theta))} = e^{i\angle(P_2(y_1))}$ but $y_1$ now belongs to a region where $\vert P_2 \vert > \vert P_1 \vert$. Hence somewhere in between there must be some $\tilde\theta$ for which $\vert P_1(y_1;\tilde\theta) \vert = \vert P_2(y_1) \vert$. This precisely corresponds to the point where the plot on Fig.~\ref{fig:lint} crosses $\vert z_* \vert = \nu = 0.8000$, and suffices to prove Theorem 3. We pursue this example towards further insight on bad or good $\theta$. Moving $e^{i\theta}$ beyond $y_1$ finally changes $\angle(P_1(y_1;\theta))$ by $\pi$ and the points along the unit circle where $e^{i\angle(P_1)} = e^{i\angle(P_2)}$ start moving; for a large interval of $\theta$ values, the points on $\gamma_0$ where $e^{i\angle(P_1)} = e^{i\angle(P_2)}$ also satisfy $\vert P_2 \vert = \vert P_1 \vert$, as on the third graph, hence the polynomial is worse than the optimal one with $M=2$. As $\theta$ approaches $\pi$ a behavior symmetric to the neighborhood of $\theta=0$ takes place (not shown).

\section{Conclusion}\label{sec:persp}

We have proved that the maximal achievable acceleration of a linear iterative map by adding $M-1$ memory slots to each state vector component, is already achieved with $M-1=1$. More precisely, this holds for the performance guarantee of a constant map about which we only know bounds $[\lm,\lp] \subset (0,+\infty)$ on the nonzero eigenvalues of the self-adjoint state update matrix $A$. The fact that the nonzero eigenvalues of $A$ can take \emph{any} values in $[\lm,\lp]$ is important for the property to hold. Better accelerations can be devised if more is known about $A$, see e.g.~Section \ref{sec:example}. 

A direct extension of our result would allow each subsystem to follow general linear dynamics; thinking of the consensus application with rational input-output transfer function at each node, this would mean characteristic polynomials of the form
$$(z-1) P(z) + \alpha \lambda_k Q(z)$$
with $P(z)$ and $Q(z)$ two polynomials. The memory slot setting \eqref{eq:AxPlus} is restricted to $Q(z) = z^{M-1}$, i.e.~the consensus algorithm uses only the latest information coming from the network. It is not clear if convergence can speed up with other $Q(z)$.

The basic consensus algorithm \eqref{eq:base} has a proven robustness to network incidents\footnote{Stability under switching networks holds: for any directed graphs, as long as $I-\alpha L$ has all components positive, thanks to the common Lyapunov function $V = \max_i x_i - \min_i x_i$~\cite{TsitsiklisThesis}; for undirected graphs, as long as $I-\alpha L$ has all eigenvalues inside the unit circle, thanks to the common Lyapunov function $\sum_i (x_i)^2$. The conditions on $I-\alpha L$ always hold if we start from a stable nominal network and ``packet drops'' occur, i.e.~failures can at times modify $L$ by setting to zero the gain associated to some links.}.
For the case \eqref{eq:AxPlus} with optimally tuned memory slot, one can construct examples where packet drops lead to instability. Thus the benefit of more memory slots might have to be reevaluated towards robustness.



\section*{Acknowledgment}

This research has been carried out within the Interuniversity Attraction Poles network DYSCO. The author wants to thank F.Ticozzi and S.Zampieri for stimulating discussions.

\end{document}